Волянский Р.С., Садовой А.В.

Днепродзержинский государственный технический университет

## КОНСТРУИРОВАНИЕ ФУНКЦИОНАЛОВ КАЧЕСТВА ДЛЯ СИСТЕМ УПРАВЛЕНИЯ С ПОКАЗАТЕЛЬНОЙ АКТИВАЦИОННОЙ ФУНКЦИЕЙ

**Введение.** Существующие системы управления электромеханическими объектами обеспечивают экспоненциально протекающие переходных процессов и гарантируют требуемый запас устойчивости [1]. Такой характер движения объекта управления в замкнутой системе объясняется в первую очередь линейностью используемых управляющих воздействий. Аналогичная ситуация наблюдается в системах разрывного управления, которые в скользящем режиме эквивалентны линейным с неограниченным коэффициентом усиления, а при срыве скользящего режима работают как разомкнутые с максимальным управляющим воздействием. За счет форсировки управляющего воздействия системы разрывного управления обладают высоким быстродействием и точностью. Однако процесс управления в таких системах сопровождается подачей на объект в каждый момент времени максимального управляющего воздействия, что не всегда является энергетически целесообразным и приводит к значительным энергетическим затратам. Поэтому, наряду с рассмотренными ранее системами в последнее время появились разработки, позволяющие форсировать протекание переходных процессов без повышения энергии управления в установившихся режимах. В первую очередь к таким разработкам следует отнести системы управления, реализующие скользящие режимы высоких порядков [2], которые возникают при использовании иррациональных активационных функций и реализуются управлениями вида

$$U = -|S|^\alpha \operatorname{sign}(S), \qquad (1)$$

где $S$ – линия равновесного состояния регулятора; $\alpha$ - показатель степени, $\alpha \in [0,1]$.

Выполненные ранее исследования [3] показали, что существенное влияние на быстродействие систем с управлениями (1) оказывает показатель степени $\alpha$, с уменьшением которого повышается не только быстродействие, но и колебательность системы. Это позволяет форсировать систему при пуске, однако негативно сказывается на установившемся режиме.

Данный факт позволил сформулировать гипотезу об изменении показателя степени $\alpha$ в процессе движения системы управления, в соответствии с которой замкнутая система при больших рассогласованиях работает как релейная с $\alpha = 0$, а по мере снижения ошибки управления увеличивается показатель $\alpha$, что устраняет автоколебания и существенно ограничивает энергию управления. Таким образом в установившемся режиме система может работать в соответствии с линейным алгоритмом управления, обеспечивая минимальное энергопотребление, а при возрастании отклонения управляемой величины от заданного значения автоматически снижается показатель $\alpha$ и за счет повышенного энергопотребления происходит форсированный переход в новую точку фазового пространства. Таким образом, за счет плавного изменения активационной функции организуется система оптимального управления с переменной структурой.

Проверке данной гипотезы в настоящее время препятствует отсутствие функционалов качества, минимизация которых осуществляется управляющим воздействием вида

$$U = -f_1(S)^{f_2(S)}. \qquad (2)$$

Управление (2) является сложной функцией от $S$ и определение интегральных функционалов качества, которые оно минимизирует на траекториях управляемого движения, является достаточно сложной задачей. Поэтому прежде чем переходить к решению этой задачи, целесообразно рассмотреть определение целей управления, достижение которых обеспечивают управляющие воздействия с показательной активационной функцией

$$U = -C^{|f_2(S)|}\text{sign}(S), \qquad (4)$$

где C – отличная от 0 и ±1 константа, которая определяет амплитуду управляющего воздействия.

**Постановка задач исследования.** Целью настоящей работы является определение интегрального функционала качества, минимизация которого осуществляется оптимальным управлением (4).

**Результаты исследования.**

Искомый функционал качества представим в виде интеграла двух слагаемых, первое из которых $F(S)$ определяет асимптотичность траекторий движения системы управления, а второе $G(U)$ – расход энергии управления

$$I = \int_0^\infty [F(S) + G(U)]dt. \qquad (5)$$

Для определения слагаемых функционала (5) в соответствии с [4] запишем сопряженные функции

$$f(S) = -g^{-1}(U) = C^{|f_2(S)|}\text{sign}(S);$$
$$g(U) = -f^{-1}(S) = \log_C|U|\text{sign}(U). \qquad (6)$$

Расход энергии на управление $G(U)$ может быть найден путем интегрирования последнего выражения системы (6) по управляющему воздействию

$$G(U) = \int g(U)dU = \int \log_C|U|\text{sign}(U)dU = \frac{1}{\ln C}\int \ln|U|\text{sign}(U)dU =$$
$$= \frac{1}{\ln C}|U|(\ln|U|-1) = C_1|U|\ln|U| - C_1|U|, \qquad (7)$$

где $C_1 = \frac{1}{\ln C}$.

Графическое представление функций $g(U)$ и $G(U)$ для случая $C = 2$ приведено на рис.1.

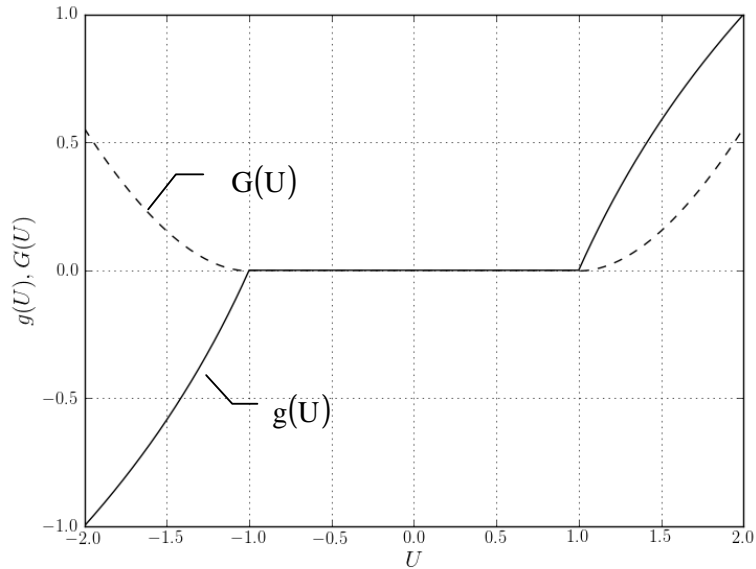

Рис.1 Внешний вид функций $g(U)$ и $G(U)$.

Анализ выражения (7) и графиков, приведенных на рис.1, показывает, что функция $G(U)$ является знакопостоянной функцией, имеющей экстремум, соответственно интеграл от нее

$$I_U = \int_0^\infty G(U)dt \qquad (8)$$

имеет экстремали и может использоваться при построении оптимальных систем.

Перейдем теперь к рассмотрению составляющей $F(S)$ функционала (5). Эта составляющая в самом общем случае определяется выражением [4]

$$F(S) = S \cdot f(S) - \int S \frac{\partial f(S)}{\partial S} dS, \qquad (9)$$

которое с учетом функции (6) примет вид

$$F(S) = S \cdot C^{|f_2(S)|} \text{sign}(S) - \int S \frac{\partial C^{|f_2(S)|}}{\partial |f_2(S)|} \frac{\partial f_2(S)}{\partial S} dS = \\
= S \cdot C^{|f_2(S)|} \text{sign}(S) - \ln C \int S C^{|f_2(S)|} \text{sign}(S) \frac{\partial f_2(S)}{\partial S} dS. \qquad (10)$$

Вычисление интеграла, образующего второе слагаемое выражения (10) в общем виде при неопределенной функции $f_2(S)$ затруднено. Поэтому имеет смысл рассмотреть определение значения выражения (10) для ряда случаев при различных функциях $f_2(S)$.

### 1. Функция $f_2(S) = S$

В этом случае выражение (10) принимает вид

$$F(S) = S \cdot C^{|S|} \text{sign}(S) - \ln C \int S C^{|S|} \text{sign}(S) \frac{\partial S}{\partial S} dS = \\
= S \cdot C^{|S|} \text{sign}(S) - \ln C \int S C^{|S|} \text{sign}(S) dS = \\
= |S| \cdot C^{|S|} - \ln C \int |S| C^{|S|} dS = |S| \cdot C^{|S|} - \ln C \frac{C^{|S|}(|S|\ln C - 1)}{(\ln C)^2} = \\
= |S| \cdot C^{|S|} - \frac{\ln C |S| C^{|S|}}{\ln C} + \frac{C^{|S|}}{\ln C} = \frac{C^{|S|}}{\ln C} = C_1 C^{|S|}. \qquad (11)$$

Очевидно, что выражение $C_1 C^{|S|}$ имеет экстремум при $S = 0$.

Подстановка значений выражений (11) и (7) в функционал (5) позволяет получить следующий функционал качества

$$I = \int_0^\infty \left[ C_1 C^{|S|} + (C_1|U|\ln|U| - C_1|U|) \right] dt = C_1 \int_0^\infty \left[ C^{|S|} + (|U|\ln|U| - |U|) \right] dt, \qquad (12)$$

который минимизируется оптимальным управлением

$$U = -C^{|S|} \text{sign}(S). \qquad (13)$$

### 2. Функция $f_2(S) = \sqrt{|S|}$

Для этой функции первое слагаемое интегранта функционала (5) примет вид

$$F(S) = S \cdot C^{\sqrt{|S|}} \text{sign}(S) - \ln C \int S C^{\sqrt{|S|}} \text{sign}(S) \frac{\partial \sqrt{|S|}}{\partial S} dS = \qquad (14)$$

$$= S \cdot C^{\sqrt{|S|}} \text{sign}(S) - \ln C \int \frac{|S|C^{\sqrt{|S|}}}{2\sqrt{|S|}} dS = |S| \cdot C^{\sqrt{|S|}} - \frac{\ln C}{2} \int \sqrt{|S|} C^{\sqrt{|S|}} dS =$$

$$= |S| \cdot C^{\sqrt{|S|}} - \ln C \frac{C^{\sqrt{|S|}}\left(|S|(\ln C)^2 - 2\sqrt{|S|}\ln C + 2\right)}{(\ln C)^3} =$$

$$= |S| \cdot C^{\sqrt{|S|}} - \ln C \frac{|S| \cdot C^{\sqrt{|S|}}(\ln C)^2}{(\ln C)^3} + \ln C \frac{2\sqrt{|S|} \cdot C^{\sqrt{|S|}} \ln C}{(\ln C)^3} - \ln C \frac{2C^{\sqrt{|S|}}}{(\ln C)^3} =$$

$$= \frac{2\sqrt{|S|} \cdot C^{\sqrt{|S|}}}{\ln C} - \frac{2C^{\sqrt{|S|}}}{(\ln C)^2} = 2\frac{C^{\sqrt{|S|}}}{\ln C}\left(\sqrt{|S|} - \frac{1}{\ln C}\right) = 2C_1 C^{\sqrt{|S|}}\left(\sqrt{|S|} - C_1\right)$$

(14)

Подставив в функционал (5) значения выражений (14) и (7), получим следующий интегральный функционал

$$I = \int_0^\infty \left[2C_1 C^{\sqrt{|S|}}\left(\sqrt{|S|} - C_1\right) + \left(C_1|U|\ln|U| - C_1|U|\right)\right] dt =$$

$$= C_1 \int_0^\infty \left[2C^{\sqrt{|S|}}\left(\sqrt{|S|} - C_1\right) + \left(|U|\ln|U| - |U|\right)\right] dt,$$

(15)

который минимизируется оптимальным управлением

$$U = -C^{\sqrt{|S|}} \text{sign}(S). \tag{16}$$

В рассмотренных примерах интеграл, образующий второе слагаемое выражения (10), определялся через элементарные функции. Однако в ряде случаев найти значение этого интеграла в классе элементарных функций не представляется возможным.

Проиллюстрируем это следующими примерами.

### 3. Функция $f_2(S) = \frac{1}{|S|}$.

Для гиперболической функции первое слагаемое интегранта функционала (5) примет вид

$$F(S) = S \cdot C^{1/|S|} \text{sign}(S) - \ln C \int S C^{1/|S|} \text{sign}(S) \frac{\partial(1/|S|)}{\partial S} dS =$$

$$= |S| \cdot C^{1/|S|} + \ln C \int |S| \frac{C^{1/|S|}}{|S|^2} dS = |S| \cdot C^{1/|S|} + \ln C \int \frac{C^{1/|S|}}{|S|} dS.$$

(17)

Найти значение интеграла, образующего второе слагаемое выражения (17), в элементарных функциях невозможно, однако значение $F(S)$ может быть определено через экспоненциальный интеграл [5]

$$\text{Ei}(a, z) = \int_1^\infty e^{-kz} k^{-a} dk \tag{18}$$

следующим образом:

$$F(S) = |S| \cdot C^{1/|S|} + \ln C \cdot \text{Ei}\left(1, -\frac{\ln C}{|S|}\right) = |S| \cdot C^{1/|S|} + \frac{1}{C_1} \text{Ei}\left(1, -\frac{1}{C_1|S|}\right). \tag{19}$$

Тогда с учетом выражений (19) и (7) искомый функционал примет вид

$$I = \int_0^\infty \left[ |S| \cdot C^{1/|S|} + \frac{1}{C_1} \mathrm{Ei}\left(1, -\frac{1}{C_1|S|}\right) + \left(C_1|U|\ln|U| - C_1|U|\right) \right] dt. \quad (20)$$

Функционал (20) минимизируется оптимальным управлением

$$U = -C^{1/|S|} \mathrm{sign}(S). \quad (21)$$

## 4. Функция $f_2(S) = |S|^\alpha$

Найдем функционал качества, который минимизируется оптимальным управлением

$$U = -C^{|S|^\alpha} \mathrm{sign}(S), \alpha \neq 0. \quad (22)$$

Использование при вычислении второго слагаемого выражения (10) верхней неполной гамма-функции [5]

$$\Gamma(s, x) = \int_x^\infty t^{s-1} e^{-t} dt, \quad (23)$$

которая связана с интегралом (18) соотношением

$$\mathrm{Ei}(a, z) = z^{a-1} \Gamma(1-a, z) \quad (24)$$

позволяет обобщить полученные результаты на случай произвольного показателя степени $\alpha$ функции

$$f_2(S) = |S|^\alpha. \quad (25)$$

Тогда выражение (10) примет вид

$$F(S) = S \cdot C^{|S|^\alpha} \mathrm{sign}(S) - \ln C \int SC^{|S|^\alpha} \mathrm{sign}(S) \frac{\partial |S|^\alpha}{\partial S} dS =$$
$$= |S| \cdot C^{|S|^\alpha} - \alpha \ln C \int |S| C^{|S|^\alpha} |S|^{\alpha-1} dS = |S| \cdot C^{1/|S|} - \alpha \ln C \int C^{|S|^\alpha} |S|^\alpha dS = \quad (26)$$
$$= |S| \cdot C^{|S|^\alpha} + \alpha \ln^{(\alpha-1)/\alpha} C \frac{\Gamma\left(\frac{\alpha+1}{\alpha}, -\ln C |S|^\alpha\right)}{\sqrt[\alpha]{-\alpha}}.$$

Подставив выражения (26) и (7) в функционал (5), получим

$$I = \int_0^\infty \left[ |S| \cdot C^{|S|^\alpha} + \alpha \ln^{(\alpha-1)/\alpha} C \frac{\Gamma\left(\frac{\alpha+1}{\alpha}, -\ln C |S|^\alpha\right)}{\sqrt[\alpha]{-\alpha}} + \left(C_1|U|\ln|U| - C_1|U|\right) \right] dt. \quad (27)$$

Дальнейшее усложнение функции $f_2(S)$, в частности переход к аддитивным формам

$$f_2(S) = \sum_{i=1}^m w_i f_i(S), \quad (28)$$

где m - число слагаемых функции $f_2(S)$, приводит к необходимости использования и других неэлементарных функций.

В частности использование функции ошибок [5]

$$\mathrm{erf}(x) = \frac{2}{\sqrt{\pi}} \int_0^x e^{-t^2} dt \quad (29)$$

позволило определить функционал, который минимизируется управлением

$$U = -C^{w_1\sqrt{|S|} + w_2|S|} \mathrm{sign}(S), \quad (30)$$

однако вследствие громоздкости полученного интеграла он здесь не приводится.

**Выводы.** Приведенные математические выкладки позволяют сделать вывод о том, что функционалы качества, которые минимизируются оптимальными управлениями с показательной активационной функцией, состоят из двух слагаемых, определяющих динамику объекта управления и энергию управления, необходимую для осуществления движения объекта по заданным траекториям возмущенного движения.

Необходимо отметить, что составляющая $G(U)$, которая учитывает расход энергии управления, однозначно зависит от управляющего воздействия и при $|U| \leq 1$ имеет экстремум, наличие которого позволяет математически обосновать возможность снижения затрат энергии управления в оптимальных системах, реализующих управляющие воздействия (4).

Составляющая $F(S)$ функционала качества, определяющая траекторию возмущенного движения объекта управления, в простейших случаях может быть выражена через элементарные функции, причем в этом случае по ее виду может быть определен класс оптимального управления. В более сложных случаях эта составляющая не может быть определена среди элементарных функций и выражается через ряд неэлементарных функций. Использование некоторых из них позволяет значительно расширить класс минимизируемых функционалов, а значит и увеличить область используемых управляющих воздействий.

Сложность приведенного математического аппарата создает предпосылки к использованию в инженерных расчетах численных методов для определения составляющих функционала (5) и алгоритмов оптимального управления (4).

Литература

1. Садовой А.В. и др. Системы оптимального управления прецизионными электроприводами/ Садовой А.В., Сухинин Б.В., Сохина Ю.В. - К.: ИСИМО, 1996. – 298с.

2. Емельянов С.В., Коровин С.К. Новые типы обратной связи: Управление при неопределенности. — М.: Наука. Физматлит, 1997, 352 с.

3. Волянский Р.С. Синтез оптимальной системы управления с иррациональной активационной функцией/ Р.С.Волянский, А.В.Садовой // Вестник НТУ «ХПИ» «Проблемы автоматизированного электропривода» (Теория и практика) вып.28. – 2010.- с.49-51

4. Волянский Р.С. Решение обратной задачи аналитического конструирования регуляторов для электромеханической системы с обобщенной активационной функцией/ Р.С.Волянский, К.А Калюжный // Сборник докладов научной конференции Достижение молодых ученых в развитии иновационных процесссов в экономике, науке и образовании» .ч.1 . Брянск: изд-во БГТУ, 2011.- С.139-140.

5. Ямке Е. и др. Специальные функции. Формулы, графики, таблицы/ Ямке Е., Эмде Ф., Леш Ф.-М.:Наука, 1964.-344с.